\documentclass[12pt]{amsart}
\usepackage{amssymb,amsfonts,amsmath,amsopn,amstext,amscd,latexsym,amsthm,enumerate,mathrsfs,bbm,color}

\usepackage[margin=36mm]{geometry}

\usepackage[all,cmtip]{xy}
\usepackage[shortlabels]{enumitem}

\newtheorem*{thm}{Theorem}

\numberwithin{equation}{section}

\begin{document}

\title[It\^{o}'s theorem and monomial Brauer characters]
{It\^{o}'s theorem and monomial Brauer characters}

\author[Xiaoyou Chen] {Xiaoyou Chen}
\address{College of Science, Henan University of Technology, Zhengzhou 450001, China}
\email{cxymathematics@hotmail.com}

\author[M. L. Lewis]{Mark L. Lewis}
\address{Department of Mathematical Sciences, Kent State University, Kent, OH 44242, USA}
\email{lewis@math.kent.edu}


\date{\today}


\begin{abstract}
Let $G$ be a finite solvable group, and let $p$ be a prime.  In this note, we prove that $p$ does not divide $\varphi(1)$ for every irreducible monomial $p$-Brauer character $\varphi$ of $G$ if and only if $G$ has a normal Sylow $p$-subgroup.
\end{abstract}

\maketitle



Throughout this paper, all groups are finite.  In \cite{LuPang2016}, Pang and Lu show that the properties of solvable groups coming from the degrees of the irreducible characters can be determined using only the degrees of the monomial irreducible characters.  In other words, for solvable groups, the monomial irreducible characters are plentiful enough to be used in the place of the irreducible characters.  

In particular, Pang and Lu prove in Theorem 1.3 of \cite{LuPang2016} if $G$ is solvable and $p$ is a prime, then $p$ does not divide $\chi(1)$ for every monomial character $\chi \in {\rm Irr} (G)$ if and only if $G$ has a normal Sylow $p$-subgroup.  Hence, Pang and Lu are able to generalize the normal Sylow subgroup portion of It\^o's theorem (Corollary 12.34 of \cite{Isaacs1}).  Note that the degrees of irreducible monomial characters for the group ${\rm SL}(2, 3)$ are $1$ and $3$, and the normal Sylow $2$-subgroup of ${\rm SL}(2, 3)$ is not abelian, so they are not able to recover the full strength of It\^o's theorem in this situation.

It\^o also proved his theorem for Brauer characters of $p$-solvable groups.  In this case also, only the normality of Sylow subgroup is recovered.  See Theorem 13.1 (b) and (c) of \cite{MaWo}.  Our main theorem is to prove for solvable groups that we only need the monomial Brauer characters to prove this result.

\begin{thm}\label{Theorem}
Let $G$ be a solvable group and let $p$ be a prime.  Then $G$ has a normal Sylow $p$-subgroup if and only if $p$ does not divide $\varphi(1)$ for every monomial Brauer character $\varphi\in {\rm IBr} (G)$.
\end{thm}

Observe that the hypothesis $G$ being solvable cannot be entirely dropped.  For example, let $G=S_{5}$ be the symmetric group of degree $5$, and $p=2$.  It is easy to see that $S_{5}$ has no subgroup of order $30$, so none of the nonlinear irreducible $2$-Brauer characters of $G$ are monomial since they have degree $4$.  Obviously, $G$ has no normal Sylow $2$-subgroup.  At this time, we have not determined if we can weaken solvable hypothesis to $p$-solvable. Our proof is motivated by the proof of Theorem 1.3 in \cite{LuPang2016}.




\begin{proof} [Proof of Theorem]
If a Sylow $p$-subgroup $P$ is normal in $G$, then $P = {\bf O}_{p}(G)$ is contained in the kernel of every irreducible $p$-Brauer character of $G$.  Thus, ${\rm IBr} (G) = {\rm IBr} (G/P) = {\rm Irr} (G/P)$.  Since $G/P$ is a $p'$-group, we have for every Brauer character $\varphi \in {\rm IBr} (G)$ that $p$ does not divide $\varphi(1)$.

Conversely, suppose that $p$ does not divide $\varphi (1)$ for every monomial Brauer character $\varphi \in {\rm IBr} (G)$.  Let $N$ be a minimal normal subgroup of $G$, and let $P$ be a Sylow $p$-subgroup of $G$.  By induction $G/N$ has a normal Sylow $p$-subgroup $PN/N$, and so, $PN$ is normal in $G$.  If $N$ is a $p$-group, then $PN = P$.  Thus, $P$ is normal in $G$ as desired.  Thus, we may assume that $N$ is an elementary abelian $q$-group for some prime $q \neq p$.  By the Frattini argument it follows that  $G = NP {\bf N}_{G} (P) = N {\bf N}_{G} (P).$  Since $N \cap {\bf N}_{G}(P)$ is normal in ${\bf N}_{G}(P)$ and $N$ is abelian, $N\cap {\bf N}_{G}(P)$ will be normal in $N {\bf N}_G (P) = G$.  The minimality of $N$ implies that either $N \le {\bf N}_G (P)$ or $N \cap {\bf N}_{G} (P) = 1$.  If $N \le {\bf N}_G (P)$, then $G = {\bf N}_G (P)$ and $P$ is normal in $G$ as desired.

We assume that $N \cap {\bf N}_G (P) = 1$.  Let $1_N \ne \lambda \in {\rm IBr} (N) = {\rm Irr} (N)$, and take $T$ to be the inertia group of $\lambda$ in $G$.  Since $N$ is complemented in $G$, it follows that $N$ is complemented in $T$.  Using Problem 6.18 of \cite{Isaacs1}, we see that $\lambda$ extends to $\nu \in {\rm Irr} (T)$.  Taking $\mu$ to be the restriction of $\nu$ to the $p$-regular elements of $T$, we see that $\mu\in {\rm IBr} (T)$ and $\mu_{N} = \lambda$.  Applying the Clifford correspondence for Brauer characters, \cite[Theorem 8.9]{Navarro1998}, we have $\varphi = \mu^{G} \in {\rm IBr} (G)$.  This implies that $\varphi$ is monomial with degree $|G: T|$.  By hypothesis, $p$ does not divide $\varphi(1) = |G: T|$.  It follows that $T$ contains some Sylow $p$-subgroup of $G$ and without loss of generality, we may assume that $P\leq T$.  Now, for all elements $x\in P$ and $n\in N$, we have $\lambda (n) = \lambda^{x} (n) = \lambda (xnx^{-1})$.  Since $\lambda$ is linear, we obtain $\lambda(xnx^{-1}n^{-1})=1$. Because $\lambda$ is arbitrary, it follows that $[P,N] \le \bigcap_{\lambda\in {\rm IBr}(N)} \ker\lambda = 1.$  Therefore, $N$ normalizes $P$.  This implies that $P$ is a characteristic subgroup of $NP$, and thus, $P$ is a normal subgroup of $G$.
\end{proof}


\section*{Acknowledgments}
The first author thanks the support of China Scholarship Council,
Department of Mathematical Sciences of Kent State University for its hospitality,
Funds of Henan University of Technology (2014JCYJ14, 2016JJSB074, 26510009),
Project of Department of Education of Henan Province (17A110004),
Projects of Zheng-zhou Municipal Bureau of Science and Technology (20150249, 20140970),
and the NSFC (11571129).



\begin{thebibliography}{99}


\bibitem{Isaacs1} I. M. Isaacs, {\it Character Theory of Finite Groups}, Academic Press, New York, 1976.


\bibitem{MaWo} O. Manz, T. R.  Wolf, \emph{Representations of solvable groups}, Cambridge University Press, Cambridge, 1993.


\bibitem{Navarro1998} G. Navarro, {\it Characters and blocks of finite groups}, Cambridge University Press, Cambridge, 1998.


\bibitem{LuPang2016} L. Pang and J. Lu, \textrm{Finite groups and degrees of irreducible monomial characters},
{\it J. Algebra Appl.} {\bf 15} (2016), 1650073 (4 pages).
\end{thebibliography}
\end{document}